\documentclass[preprint]{elsarticle}
\usepackage{tikz-cd}
\usepackage{cmap} % for russuan cope paste
           \usepackage[T2A]{fontenc}   % Or any other encoding.
           \usepackage[utf8]{inputenc}
           \usepackage[english]{babel}
\usepackage{amsmath}
\usepackage{amsfonts}
\usepackage{amssymb}
\usepackage{amsthm}
\usepackage[pdftex,unicode]{hyperref}
\usepackage{indentfirst} % включить отступ у первого абзаца
\usepackage{amscd}
\usepackage[all,cmtip]{xy}
\usepackage{mathtools}
\usepackage{comment}
\usepackage{tikz-cd}
\usetikzlibrary{babel}

  \def\R{\mathbb R}

\def\om{\omega}

\def\be{\beta}

\def\al{\alpha}

\def\ph{\varphi}

\def\nfs/{NFS}
\def\cdp/{CDP}
\def\cdpz/{CDP${}_0$}

\def\sq#1#2{(#1)_{#2}}
\def\sqn#1{\sq{#1}{n\in\om}}
\def\sqnn#1{\sqn{#1_n}}

\def\cl#1{\overline{#1}}

\def\sset#1{\{#1\}}

\def\set#1{\bbset#1\eeset}
\def\bbset#1:#2\eeset{\{#1\,:\,#2\}}

\def\bbsett#1:#2\eesett{\{#1\,:\,\text{#2}\}}

\def\ibbset#1:#2\ieeset{(#1)_{#2}}

\def\cP{{\mathcal P}}

\newcommand\restrA[2]{{% we make the whole thing an ordinary symbol
  \left.\kern-\nulldelimiterspace % automatically resize the bar with \right
  #1 % the function
  \vphantom{\big|} % pretend it's a little taller at normal size
  \right|_{#2} % this is the delimiter
  }}

\newcommand\restrB[2]{\ensuremath{\left.#1\right|_{#2}}}

\def\restr#1#2{\restrB{#1}{#2}}

\def\pwr#1_#2{#1^{[#2]}}

\def\term#1{{\it #1}}

\def\dddm#1(#2){N_{#1}(#2)}
\def\dddb#1(#2){B_{#1}(#2)}

\def\wh#1{\widehat{#1}}
\def\et(#1){ (#1)}

\def\bitem#1,#2.{ $#2\nrightarrow #1$:\ }
\newtheorem{proposition}{Proposition}
\newtheorem{theorem}{Theorem}

\newtheorem*{lemma*}{Lemma}
\newtheorem{cor}{Corollary}
\newtheorem{problem}{Problem}

\theoremstyle{definition}

\def\oo#1/{$O_{#1}$}

\def\gd/{$G_\delta$}

\def\rarr{\Rightarrow}
\def\rarr{\rightarrow}

\def\wh#1{\widehat{#1}}

\def\ii{{\mathfrak{i}}}

\begin{document}

\begin{frontmatter}

\title{All countable subsets of pseudocompact quasitopological Korovin groups are 
closed, discrete and $C^\ast$-embedded}

\author{Evgenii Reznichenko}\corref{cor1}
\ead{erezn@inbox.ru}
\address{Department of General Topology and Geometry, Mechanics and Mathematics Faculty, M.~V.~Lomonosov Moscow State University, Leninskie Gory 1, Moscow, 199991 Russia}

\author{Mikhail Tkachenko\fnref{fn2}}
\ead{mich@xanum.uam.mx}
\address{Departamento de Matem\'{a}ticas, Universidad Aut\'onoma Metropolitana, Av. San Rafael Atlixco 186, Col. Vicentina, C.P. 09340, Del. Iztapalapa, Mexico City, Mexico}

\cortext[cor1]{Corresponding author}

\fntext[fn2]{The second listed author was supported by grant number CAR-64356 of the 
Program \lq\lq{Ciencia de Frontera 2019\rq\rq} of the CONACyT, Mexico.}
%% or include affiliations in footnotes:

\date{May 15, 2023}

\begin{abstract}
We show that all countable subsets of any pseudocompact quasitopological group 
in the form of a Korovin orbit are closed, discrete, and $C^\ast$-embedded. Consequently, 
any infinite pseudocompact Korovin orbit is not homeomorphic to a topological group. 
Moreover, infinite pseudocompact Korovin orbits are not homeomorphic to any 
Mal'tsev space.
\end{abstract}

\begin{keyword}
Pseudocompact groups
\sep
Quasitopological groups
\sep
Korovin orbit
\sep
Mal'tsev space
\end{keyword}
\end{frontmatter}

\section{Introduction}
\label{sec:intro}
	A group with a topology is called \term{semitopological} if multiplication in the group is separately continuous.
	Ellis \cite{Ellis1957-2} proved that a locally compact semitopological group is a topological group, that is, multiplication and inverse are continuous. In \cite{Reznichenko1994}, it is proved that a completely regular countably compact semitopological group is a topological group.

	A semitopological group is called \term{quasitopological} if inversion in the group 
is continuous. In \cite{Korovin1992}, A.~Korovin presented a method for the construction 
of pseudocompact semitopological groups that are neither topological groups nor 
homeomorphic to any topological group (the latter is shown in our Corollary~\ref{Cor:tg}). 
An analysis of Korovin's method led to the notions of \emph{Korovin map} and 
\emph{Korovin orbit} \cite{ArhangelskiiHusek2001,HernandezTkachenko2006} 
(see also Theorem~2.4.13 in \cite{at2009}). Korovin orbits are also used to construct pseudocompact quasitopological groups with various additional properties \cite{ArhangelskiiHusek2001,HernandezTkachenko2006,Batikova2009,TangLinXuan2021}.
	
	Our main result states that all countable subsets of a pseudocompact 
Korovin orbit are closed, discrete, and $C^\ast$-embedded (Theorem~\ref{t:korovin:1}).
Consequently, pseudocompact Korovin's orbits are not homeomorphic to topological 
groups or even paratopological groups. Moreover, pseudocompact Korovin's orbits are 
not homeomorphic to any Mal'tsev space (Corollary~\ref{c:maltcev:1}).

	\paragraph{Around Grothendieck's theorem}
We say that a Tychonoff space $X$ is a \term{Korovin space} if for any continuous 
mapping $f\colon X\to C_p(X)$, the image $f(X)$ has compact closure. Here 
$C_p(X)$ denotes the space of continuous real-valued functions on $X$ with the 
topology of pointwise convergence. Since $C_p(X)$ contains a closed copy of 
the real line, any Korovin space is pseudocompact. In \cite{Reznichenko1994}, 
it was proved that if a semitopological group $G$ is a Korovin space, then it is 
a topological group.
	
	A Grothendieck's theorem in \cite{grot1952} states that the closure of a countably compact subspace of $C_p(X)$, for a countably compact $X$, is compact. It follows from this result that countably compact spaces are Korovin spaces, which implies the aforementioned result about countably compact semitopological groups. Numerous generalizations of Grothendieck's theorem show that the class of Korovin spaces is quite wide, for example, it includes both separable and Fréchet--Urysohn pseudocompact spaces  \cite{Korovin1992,Reznichenko1994,Arhangelskii1997}.
	
As shown in Corollary~\ref{Cor:tg},
any infinite pseudocompact Korovin orbit is not homeomorphic to a topological group.
%Hence an infinite pseudocompact Korovin orbit is not a topological group. 
According to \cite{Reznichenko1994} and the above-mentioned result, every infinite pseudocompact Korovin orbit is not a Korovin space.
	We conclude, therefore, that Grothendieck's theorem cannot be extended to pseudocompact spaces. Using D. Shakhmatov's example of an infinite pseudocompact space $X$ in which every countable subset is closed, discrete, and $C^\ast$-embedded may provide another example of this type, as pointed out by V.~Tkachuk.	
	The subspace $Y=C_p(X,[0,1])$ of $C_p(X)$ is closed in $C_p(X)$, pseudocompact, and non-compact. Theorem~\ref{t:korovin:1} shows that in the above argument, pseudocompact Korovin orbits can be used in place of Shakhmatov's space from \cite{Shakhmatov1986}.
	
	\paragraph{Mal'tsev spaces}
	In a pseudocompact topological group, any infinite subgroup is not discrete \cite{HernandezTkachenko2006}. Hence, an infinite pseudocompact Korovin orbit is not homeomorphic to a topological group. A ternary operation $M\colon X^3\to X$ is called a \term{Mal'tsev operation} if the equality
	\[
	M(y,y,x)=M(x,y,y)=x
	\]
holds for all $x,y\in X$. A space with a continuous Mal'tsev operation is called \term{Mal'tsev space}. There exists a natural Mal'tsev operation on every group $X$ defined by $M(x,y,z)=xy^{-1}z$. If $X$ is a topological group, then $M$ is continuous. If $X$ is a quasitopological group, then $M$ is separately continuous. Therefore, every pseudocompact quasitopological Korovin orbit admits a separately continuous Mal'tsev operation, but it does 
not admit a continuous Mal'tsev operation.
	
It is worth mentioning that if $X$ is a pseudocompact Mal'tsev space or a Korovin space with a separately continuous Mal'tsev operation, then the Stone--\v{C}ech compactification $\be X$ of $X$ is a Dugundji compactum \cite{ReznichenkoUspenskij1998,Reznichenko2022-2}.

\section{Korovin orbits}
\label{sec:korovin}
	Let $(G,+)$ be an Abelian group and $X$ a space with at least two elements.
	The group $G$ acts on $X^G$ as follows:
	\[
	s: G\times X^G \to X^G,\ s(g,f)(h)=f(g+h),
	\]
where $f\in X^G$ and $g,h\in G$.
	An element $f\in X^G$ is called a \term{Korovin mapping} if for any $M\subset G$ 
	with $|M|\leq\om$ and $h\in X^M$ there exists an element $g\in G$ such that $h=\restr{s(g,f)}M$.
	
	\begin{theorem}[\cite{Korovin1992} and for $|G|=2^\om$, 
	Theorem~2.4.13 in \cite{at2009}]\label{t:korovin:korovin}
		If $G$ is an Abelian group, $X$ is a space, and $2\leq |X|\leq |G|=|G^\om|$, 
		then there exists a Korovin mapping $f\colon G\to X$.
	\end{theorem}
	
	The subspace $G_f = \set{s(g,f): g\in G}$ of $X^G$,
	for a Korovin mapping $f$, is called a \term{Korovin orbit}. 
	Let $\pi_M\colon X^G\to X^M$ be the projection, where $M\subset G$.
	
	\begin{proposition}[\cite{Korovin1992} and Proposition~2.4.14 in \cite{at2009}]\label{p:korovin:0}
		If $G$ is an Abelian group, $X$ is a space, and $f\in X^G$ is a Korovin mapping, 
		then $\pi_M(G_f)=X^M$, for any countable set $M\subset G$.
	\end{proposition}
	
	Proposition~\ref{p:korovin:0} implies the following result.
	
	\begin{proposition}[\cite{Korovin1992,ArhangelskiiHusek2001}]\label{p:korovin:0+1}
	A Korovin orbit $G_f\subset X^G$ is pseudocompact if and only if $X^\om$ is pseudocompact. In this case, $\beta G_f=(\beta X)^G$.
	\end{proposition}
	
    Let $G_f\subset X^G$ be a Korovin orbit. The mapping $\ii\colon G\to G_f$ 
defined by $\ii(g) = s(g,f)$ for all $g\in G$ is bijective \cite{ArhangelskiiHusek2001}. 
We consider the algebraic structure on $G_f$ induced by $G$. This means that 
we define the sum in $G_f$ be letting $\ii(g)+\ii(h)=\ii(g+h)$ for all $g,h\in G$. 
If $e$ is the identity element of $G$, then $\ii(e)=s(e,f)$ is the identity of the group 
$G_f$. Also, $G_f$ is a semitopological group provided it carries the topology inherited 
from $X^G$ and if $G$ is a Boolean group, then $G_f$ is a quasitopological group.
	
	\begin{proposition}\label{p:korovin:1}
		Let $G$ be an Abelian group, $X$ a space, and $f\in X^G$ a Korovin mapping.
		The function $h\circ \ii^{-1}\colon G_f\to X$ is continuous, for each $h\in G_f$.
	\end{proposition}
	
	\begin{proof}
	Given an element $h\in G_f$, there is $\al\in G$ such that $h=s(\al,f)$.\smallskip
		%\begin{lemma*}

\noindent
{\bf Claim.} \emph{$(h\circ \ii^{-1})(x)=x(\al)$ for each $x\in G_f$.}\smallskip
		%\end{lemma*}
		%\begin{proof}
		
	Indeed, take $\be\in G$ such that $x=s(\be,f)$. Then
			\begin{align*}
				(h\circ \ii^{-1})(x)&=h(\ii^{-1}(x))=h(\ii^{-1}(s(\be,f))) =h(\be)=s(\al,f)(\be)
				\\
				&=f(\al+\be)=f(\be+\al)=s(\be,f)(\al)=x(\al).
			\end{align*}
		%\end{proof}
		It follows from our Claim that $h\circ \ii^{-1}=\restr{\pi_\al}{G_f}$, where
		\[
		\pi_\al: G^X\to X,\ h\mapsto h(\al)
		\]
is the projection. Since $\pi_\al$ is continuous, the mapping $h\circ \ii^{-1}$ is continuous as well.
	\end{proof}
	
	Recall that a subset $Y$ of a space $X$ is called \term{$C$-embedded} 
(\term{$C^\ast$-embedded}) in $X$ if any (bounded) continuous real-valued function on $Y$ extends to a continuous real-valued function on $X$. Two nonempty subsets $A,B\subset X$ are called \emph{completely separated} in $X$ if there exists a continuous function $f\colon X\to \R$ such that $f(A)=\sset 0$ and $f(B) =\sset 1$. According to \cite[1.17]{gli-jer-1960}, a subset $Y$ of a space
$X$ is $C^\ast$-embedded in $X$ if and only if every two completely separated subsets of
$Y$ are completely separated in $X$. It is essential to note that the space $X$ is not assumed to satisfy any separation restriction in this fact (known as \emph{Urysohn's extension theorem}).

The following proposition is immediate from Urysohn's extension theorem.

	\begin{proposition}\label{p:korovin:2}
Let $X$ be a space and $D\subset X$. Then $D$ is discrete and $C^\ast$-embedded 
in $X$ if and only if any two nonempty disjoint subsets of $D$ are completely separated
in $X$.
	\end{proposition}

%	\begin{proof}
%The necessity of the condition is evident. Conversely, if any two disjoint subsets of $M$ are %completely separated, then $M$ is discrete and the closure $K$ of the set $M$ in the Stone--%\v{C}ech compactification $\beta X$ of $X$ is the Stone--\v{C}ech compactification $\beta M$ 
%of the discrete space $M$. Then $M$ is $C^\ast$-embedded in $\beta X$ and hence in $X$.
%	\end{proof}
	
     The following theorem is the article's main result.	In it, we do not impose any separation
     restriction on $X$.
	
	\begin{theorem}\label{t:korovin:1}
		Let $G$ be an Abelian group, $X$ be a space with $|X|\geq 2$, and 
		$f\colon G\to X$  a Korovin mapping.
		\begin{enumerate}
			\item
			If $X$ is not anti-discrete, then then every countable subset of the Korovin 
			orbit $G_f$ is closed and discrete in $G_f$.
			\item
			If there exists a non-constant continuous real-valued function on $X$, then 
			every countable subset in the Korovin orbit $G_f$ is $C^\ast$-embedded.
			\item
			If there is an unbounded continuous real-valued function on $X$, then 
			any countable subset of the Korovin orbit $G_f$ is $C$-embedded.
		\end{enumerate}
	\end{theorem}
	
	\begin{proof}
		(1)
		Let $A$ be a countable subset of $G_f$ and $g\in G_f\setminus A$.
		It suffices to check that $g\notin \cl A$. Let $U$ be a proper nonempty 
		open subset of $X$. Take elements $x\in U$ and $y\in X\setminus U$.
		Since $f$ is a Korovin mapping, there exists $h\in G_f$ such that 
		$h(\ii^{-1}(g))=x$ and $h(\ii^{-1}(A))=\sset y$. Proposition~\ref{p:korovin:1} 
		implies that the mapping $h\circ \ii^{-1}$ is continuous, so 
		\[
		V=(h\circ \ii^{-1})^{-1}(U)\ni g
		\]
		is open in $G_f$ and $A\subset G_f\setminus V$. Hence $g\notin \cl A$.
		
		(2)
	According to Proposition~\ref{p:korovin:2} it suffices to prove that any
	countable disjoint subsets $A,B$ of $G_f$ are completely separated 
	in $G_f$. Let $A'=\ii^{-1}(A)$ and $B'=\ii^{-1}(B)$.
		Let $q\colon X\to \R$ be a non-constant continuous function,
		$x_1,x_2\in X$ and $q(x_1)\neq q(x_2)$.
		Let $M=A'\cup B'$ and $h\colon M\to X$ be a function such that 
		$h(A^\prime)=\sset{x_1}$ and $h(B^\prime)=\sset{x_2}$. Since
		$f$ is a Korovin mapping, there exists $g\in G$ such that $h=\restr{s(g,f)}M$.
		Proposition~\ref{p:korovin:1} implies that the mapping $\ph=s(g,f)\circ \ii^{-1}: 
	        G_f\to X$ is continuous. Then $\ph(A)=\sset{x_1}$ and $\ph(B)=\sset{x_2}$. 		  
	        The function $\psi=q\circ \ph$ is continuous, $\psi(A)=\sset{q(x_1)}$, 
	        $\psi(B)=\sset{q(x_2)}$, and $q(x_1)\neq q(x_2)$.
		
		(3)
		Let $A=\set{a_n:n\in\om}\subset G_f$ be a countable set and 
		$v\colon A\to \R$ be a function. Let also $q\colon X\to \R$ be an unbounded 
	         continuous function. There is a sequence $\sqnn x \subset X$ such that 
	         $Q=\{q(x_n): n\in\omega\}$ is a closed discrete subset of $\R$ and 
	         $q(x_n)\neq q(x_m)$ if $n\neq m$. Since $Q$ is $C$-embedded in $\R$,
		there exists a continuous function $r\colon \R\to \R$ such that $r(q(x_n))=v(a_n)$ 
		for each $n\in\om$. Let $M=\ii^{-1}(A)$ and $h\colon M\to X$ be a function such 
		that $h(\ii^{-1}(a_n))=x_n$ for each $n\in\om$. Since $f$ is a Korovin mapping, 
		there exists $g\in G$ such that $h=\restr{s(g,f)}M$. Proposition~\ref{p:korovin:1} 
		implies that the mapping $\ph=s(g,f)\circ \ii^{-1}\colon G_f\to X$ is continuous. 
		Then the mapping $\psi= r\circ q \circ \ph$ is continuous and $v = \restr \psi A$.
	\end{proof}
	
     For the reader convenience, we present a short proof of the following well-known
     result (see e.g. \cite[Corollary~1.4.24]{at2009}).
	
	\begin{proposition}\label{p:korovin:css}
	An infinite pseudocompact topological group has a countable non-closed subset.
	\end{proposition}
	
\begin{proof}
A pseudocompact topological group is precompact \cite[Theorem~1.1]{CR}, a subgroup 
of a precompact topological group is precompact, and a discrete precompact group is finite. 
Hence, in an infinite pseudocompact topological group $G$, each infinite subgroup is not discrete, and $G$ has a countable non-discrete subspace. 
\end{proof}	
	
In fact, the above argument shows that Proposition~\ref{p:korovin:css} remains valid for 
infinite precompact topological groups. 
	
Proposition~\ref{p:korovin:css} and Theorem~\ref{t:korovin:1} imply the following assertion
announced in the abstract.

           \begin{cor}\label{Cor:tg}
      Any infinite pseudocompact Korovin orbit is not homeomorphic to a topological group.
          \end{cor}

\section{Discrete subsets of Mal'tsev spaces}
\label{sec:maltcev}
	
	Our aim in this short section is to extend Corollary~\ref{Cor:tg} to pseudocompact
	Mal'tsev spaces.
	
		A subset $Y$ of a Mal'tsev space $X$ with a Mal'tsev operation $M$ is called a \term{Mal'tsev subalgebra} if $Y$ is closed under the operation $M$, that is $M(Y^3)=Y$.
	\begin{theorem}\label{t:maltcev:1}
		Let $X$ be a pseudocompact Mal'tsev space and $Y$ be an infinite discrete Mal'tsev subalgebra of $X$. Then the closure of $Y$ in the Stone--\v{C}ech compactification $\beta X$ 
of $X$ is a metrizable compact set and $Y$ is not $C^\ast$-embedded in $X$.
	\end{theorem}
	
	\begin{proof}
		Let $M$ be a continuous Mal'tsev operation on $X$. Then $M$  
		extends to a continuous Mal'tsev operation $\wh M\colon (\beta X)^3\to\beta X$ 
		on the Stone--\v{C}ech compactification $\beta X$ of $X$ \cite{ReznichenkoUspenskij1998}. Let $K$ be the closure of $Y$ in $\beta X$. Then $K$ is a Mal'tsev subalgebra of $\beta X$ and $K$ is a Mal'tsev space with the dense discrete subspace $Y$. Compact Mal'tsev spaces are Dugundji compact spaces \cite{Uspenskii1989,ReznichenkoUspenskij1998}. Dugundji compacta are dyadic compacta, 
and every dyadic compact space with a dense set of isolated points is metrizable  (see \cite{efimov1965} or \cite[Corollary~12]{Arkhangelskii1989}). Hence, $K$ is a metrizable 
compact space. If $Y$ is $C^\ast$-embedded in $X$, then $K$ is homeomorphic to $\beta Y$ and $\beta Y$ is not metrizable, which is a contradiction. 
	\end{proof}
	Every countable subset of a Mal'tsev space $X$ is contained in a countable
        Mal'tsev subalgebra of $X$. Therefore, Theorem~\ref{t:maltcev:1}  
        implies the following fact.
	\begin{proposition}\label{p:maltcev:css}
	An infinite pseudocompact Mal'tsev space contains either a countable non-closed 
	subset or a countable subset that fails to be $C^\ast$-embedded.	
	\end{proposition}
	Combining Proposition~\ref{p:maltcev:css} and Theorem~\ref{t:korovin:1}
        we deduce the following assertion that improves upon Corollary~\ref{Cor:tg}.        
	\begin{cor}\label{c:maltcev:1}
		An infinite pseudocompact Korovin orbit is not homeomorphic to any 
		pseudocompact Mal'tsev space.
	\end{cor}
	We do not know whether every infinite pseudocompact Mal'tsev space 
	contains a countable non-closed subset (see Problem~\ref{pr:maltcev:1} 
	below). In the following proposition we establish that pseudocompact 
	Mal'tsev spaces have a weaker property.
	\begin{proposition}\label{p:maltcev:css2}
	An infinite pseudocompact Mal'tsev space contains a non-closed subset 
	of cardinality less than or equal to $\om_1$.
	\end{proposition}
	\begin{proof}
		Let $X$ be an infinite pseudocompact Mal'tsev space. If $X$ is countable,
		then it is compact and metrizable, so the conclusion of the proposition is 
		evident. We assume therefore that $|X|\geq\omega_1$. Let $Y$ be a Mal'tsev 
		subalgebra of $X$ of cardinality $\om_1$. Then $Y$ is not discrete\,---\,otherwise
                Theorem~\ref{t:maltcev:1} implies that the closure of $Y$ in the Stone--\v{C}ech 
		compactification $\beta X$ of $X$ is a metrizable compact space with an 
	        uncountable discrete subspace. Consequently, $Y$ is not discrete 
	        and the set $Y\setminus \{y\}$ is not closed for some non-isolated point
	        $y$ in $Y$.
	\end{proof}

\section{Questions}
\label{sec:questions}
In connection with Propositions~\ref{p:korovin:css}--\ref{p:maltcev:css2} 
%and~\ref{p:maltcev:css} 
it is natural to pose the following problem.
	\begin{problem}\label{pr:maltcev:1}
	Does there exist an infinite pseudocompact Mal'tsev space in which every countable   
	subset is closed and discrete?
	\end{problem}
 We mention the following properties of the pseudocompact space $X$ constructed in \cite{rezn1989} in relation with Problem~\ref{pr:maltcev:1} and Theorem~\ref{t:maltcev:1}: 
 (1) $X$ is a pseudocompact space of cardinality continuum , and each subset of $X$ whose cardinality is less than the continuum is closed and discrete; (2) $\beta X$ is homeomorphic to a Tychonoff cube of uncountable weight and is a Dugundji compactum; (3) the closure of every countable subset of $X$ in $\beta X$ is a metrizable compact space; (4) $X$ is a hereditarily metalindel\"of space 
 \cite{Pavlov2009}.
Proposition~\ref{p:maltcev:css2} implies that, assuming the negation of the Continuum Hypothesis, the space $X$ cannot be a Mal'tsev space. If we assume the Continuum 
Hypothesis, then it is not clear whether the space $X$ is Mal'tsev or not.

%\footnote{{\color{red}Is it known that the space $X$ is NOT a Mal'tsev space? In any case, we have to explain explicitly what is known about this space.}}

Considering Corollaries~\ref{Cor:tg} and~\ref{c:maltcev:1}, the following problem arises.
	\begin{problem}\label{pr:maltcev:2}
	Can an infinite Korovin orbit $G_f\subset X^G$ be homeomorphic to a topological 
	group or a Mal'tsev space? What if $X$ is a (countable) discrete space?
	\end{problem}

A space $X$ is called \term{weakly pc-Grothendieck} if the closure of every pseudocompact subspace $C_p(X)$ is compact \cite{Arhangelskii1997}. A pseudocompact weakly pc-Grothendieck space is a Korovin space.

The only instances of pseudocompact semitopological groups that are not topological 
groups are Korovin orbits. Also, Shakhmatov's example in \cite{Shakhmatov1986} is a pseudocompact space that is not weakly pc-Grothendieck because the closed 
pseudocompact subspace $Y=\{f\in C_p(X): f(X)\subset [-1,1]\}$ of $C_p(X)$ is not 
compact.

Let us call a pseudocompact space $Z$ \term{Shakhmatov space} if any countable subset of $Z$ is discrete, closed, and $C^\ast$-embedded. The space in \cite{Shakhmatov1986} and pseudocompact Korovin orbits are Shakhmatov spaces. 

\begin{problem}\label{pr:questions:2}
Let $\cP$ be one of the classes listed below.
\begin{enumerate}
\item[{\rm ($\cP_1$)}]
Pseudocompact spaces with countable extend.
\item[{\rm ($\cP_2$)}]
Pseudocompact spaces that have a dense subspace with countable extend.
\item[{\rm ($\cP_3$)}]
Pseudocompact spaces that have a dense Lindel\"of subspace.
\end{enumerate}
Find out which of the following statements are true:
\begin{enumerate}
\item
There is a Korovin orbit from $\cP$.
\item
There is a Shakhmatov space from $\cP$.
\item
There is a not weakly pc-Grothendieck space from $\cP$.
\item
There is a semitopological group from $\cP$ which is not a topological group.
\end{enumerate}
\end{problem}

It is easy to see that $(1)\rarr (2)\rarr (3)$ and $(1)\rarr (4)\rarr (3)$.

\begin{problem}\label{pr:questions:3}
Let $\cP$ be one of the classes listed below.
\begin{enumerate}
\item[{\rm ($\cP_4$)}]
Pseudocompact spaces that have a dense subspace with countable tightness.
\item[{\rm ($\cP_5$)}]
Pseudocompact spaces with a dense subspace which is $k$-space.
\end{enumerate}
Find out which of the following statements are true:
\begin{enumerate}
\item
There is a space from $\cP$ that fails to be weakly pc-Grothendieck.
\item
There is a semitopological group from $\cP$ which is not a topological group.
\end{enumerate}
\end{problem}
Note that in Problem~\ref{pr:questions:3}, (2) implies (1). We also mention that Shakhmatov's space in \cite{Shakhmatov1986} does not help to answer (1) of Problem~\ref{pr:questions:3}, while Korovin orbits do not help answering (2) of the problem. These conclusions follow from 
the fact that both Shakhmatov spaces of countable tightness and Shakhmatov $k$-spaces are discrete.

\bibliographystyle{elsarticle-num}
\bibliography{pseudocompact_groups6}
\end{document}